\documentclass[11pt]{amsart}
\usepackage{amssymb,amsmath,amsthm}
\usepackage{a4}
\usepackage{epsf}
\usepackage{epsfig}
\usepackage{url}
\begin{document}

\newcommand{\NP}{$\mathcal{N}\mathcal{P}$}
\newcommand{\newsymb}{{\mathcal P}}
\newcommand{\Pn}{{\mathcal P}^n}
\newcommand{\R}{{\mathbb R}}
\newcommand{\N}{{\mathbb N}}
\newcommand{\Q}{{\mathbb Q}}
\newcommand{\Z}{{\mathbb Z}}
\newcommand{\C}{{\mathbb C}}

\newcommand{\enorm}[1]{\Vert #1\Vert}
\newcommand{\inter}{\mathrm{int}}
\newcommand{\conv}{\mathrm{conv}}
\newcommand{\aff}{\mathrm{aff}}
\newcommand{\lin}{\mathrm{lin}}
\newcommand{\cone}{\mathrm{cone}}
\newcommand{\bd}{\mathrm{bd}}

\newcommand{\dist}{\mathrm{dist}}
\newcommand{\trans}{\intercal}
\newcommand{\diam}{\mathrm{diam}}
\newcommand{\vol}{\mathrm{vol}}
\newcommand{\LE}{\mathrm{G}}
\newcommand{\lE}{\mathrm{g}}
\newcommand{\sa}{\mathrm{a}}

\newcommand{\pp}{\mathfrak{p}}
\newcommand{\pf}{\mathfrak{f}}
\newcommand{\pg}{\mathfrak{g}}
\newcommand{\PP}{\mathfrak{P}}
\newcommand{\pl}{\mathfrak{l}}
\newcommand{\pv}{\mathfrak{v}}
\newcommand{\cl}{\mathrm{cl}}
\newcommand{\bx}{\overline{x}}

\def\ip(#1,#2){#1\cdot#2}

\newtheorem{theorem}{Theorem}[section]
\newtheorem{theorem*}{Theorem}
\newtheorem{corollary}[theorem]{Corollary}
\newtheorem{lemma}[theorem]{Lemma}
\newtheorem{remark}[theorem]{Remark}
\newtheorem{definition}[theorem]{Definition}  
\newtheorem{conjecture}{Conjecture}[section] 
\newtheorem{proposition}[theorem]{Proposition}  
\newtheorem{claim}[theorem]{Claim}
\newtheorem{problem}[theorem]{Problem}
\numberwithin{equation}{section}

\title[Notes on the roots of Ehrhart polynomials]{Notes on the
  roots of Ehrhart polynomials}   
\author{Christian Bey}
\address{Christian Bey, Universit\"at Magdeburg, IAG,
  Universit\"ats\-platz 2, D-39106 Magdeburg, Germany}
\email{christian.bey@mathematik.uni-magdeburg.de}
\author{Martin Henk}
\address{Martin Henk, Universit\"at Magdeburg, IAG,
  Universit\"ats\-platz 2, D-39106 Magdeburg, Germany}
\email{henk@math.uni-magdeburg.de}
\author{J\"org M.~Wills}
\address{J\"org M.~Wills, Universit\"at Siegen,  Mathematisches Institut,
 ENC, D-57068 Siegen, Germany}
\email{wills@mathematik.uni-siegen.de}


\begin{abstract} We determine lattice polytopes of smallest volume with a given number of interior lattice points. We show that the Ehrhart polynomials of those with one interior lattice point have largest roots with norm of order $n^2$, where $n$ is the dimension. This improves on the previously best known bound $n$ and complements a recent result of Braun \cite{braun:ehrhart_roots} where it is shown that the norm of a root of a Ehrhart polynomial is at most of order $n^2$. 

For the class of $0$-symmetric lattice polytopes we present a conjecture on the smallest volume for a given number of interior lattice points and prove the conjecture for crosspolytopes.
    
We further give  a characterisation of the roots  of the Ehrhart polyomials in the
  $3$-dimensional case and we classify for $n\leq 4$ all lattice
  polytopes whose roots of their Ehrhart polynomials have all real part
  -1/2. These polytopes belong to  the class of reflexive
  polytopes. 
\end{abstract}

\maketitle

\section{Introduction}
Let $\mathcal{P}^n$ be the set of all convex lattice $n$-polytopes in the $n$-dimensional Euclidean space $\R^n$ with respect to the standard lattice $\Z^n$, i.e., all vertices of $P\in \mathcal{P}^n$ have  integral coordinates and $\dim(P)=n$. The lattice point enumerator of a set $S\subset\R^n$ is denoted by $\LE(S)$,  i.e., $\LE(S)=\#(S\cap \Z^n)$. 

In 1962 Ehrhart \cite{ehrhart:polynomial} showed that for $k\in \N$ the lattice point enumerator $\LE(k\,P)$, $P\in \mathcal{P}^n$, is a polynomial of degree $n$ in $k$ where the coefficients $\LE_i(P)$, $0\leq i\leq n$,  depend only on $P$:
\begin{equation}
   \LE(k\,P)=\sum_{i=0}^n \LE_i(P)\,k^i.
\label{eq:ehrhart_polynomial}
\end{equation} 
Moreover in \cite{ehrhart:reciprocity} he  proved his famous ``reciprocity law''
\begin{equation}
   \LE(\inter(k\,P))=(-1)^{n} \sum_{i=0}^n \LE_i(P)\,(-k)^i,
\label{eq:reciprocitylaw}
\end{equation}
where $\inter()$ denotes the interior. Two  of the $n+1$ coefficients
$\LE_i(P)$ are obvious, namely, $\LE_0(P)=1$ and $\LE_n(P)=\vol(P)$,
where $\vol()$ denotes the volume, i.e., the $n$-dimensional Lebesgue
measure on $\R^n$.  Also the second leading coefficient admits a
simple geometric interpretation  as normalized surface area of $P$ which we present in detail in \eqref{eq:2leadterm}. 
All other coefficients $\LE_i(P)$, $1\leq i\leq n-2$, have no such
direct geometric meaning, except for special classes of polytopes
(cf., e.g., \cite{beckpixton:birkhoff, betkegritzmann:valuations, diazrobins:latticepolytope,  hibi:92, liu:lattice_face, liu:cyclic_ehrhart,  mordell:dedkind,  mustatepayne:bettyehrhart, pommersheim:tetrahedra}).

A sometimes more convenient representation of $\LE(k\,P)$ is given by
a change from the monomial basis $\{x^i: i=0,\dots,n\}$ to the basis
$\{\binom{x+n-i}{n}: i=0,\dots,n\}$: 
\begin{equation}
 \LE(k\,P) = \sum_{i=0}^n \sa_i(P) \binom{k+n-i}{n}.
\label{eq:ehrhart_series}
\end{equation} 
In view of \eqref{eq:ehrhart_polynomial} and \eqref{eq:reciprocitylaw}
we get 
\begin{equation}
\begin{split}
  &\sa_0(P)=1,\quad\sa_1(P)=\LE(P)-(n+1),\quad \sa_n(P)=\LE(\inter(P)),\\
  &\sa_0(P)+\sa_1(P)+\ldots+\sa_n(P)=n!\,\vol(P),
\end{split}
 \label{eq:stanley_coefficients}
\end{equation} 
and all $\sa_i(P)$ are integers. 
Due to Stanley's famous  non-negativity theorem
\cite{stanley:nonnegative} they are also 
non-negative, in contrast to the $\LE_i(P)$'s
which might be negative.

In recent years the Ehrhart polynomial was not only regarded as a polynomial for integers $k$, but as a formal polynomial of a complex variable $s\in \C$ (cf.~\cite{becketall:roots_ehrhart, rodriguez:zeros, wills:90c}). Therefore, for $P\in\Pn$ and $s\in \C$ we set   
\begin{equation*}
   \LE(s, P)= \sum_{i=0}^n \LE_i(P)\, s^i = \prod_{i=1}^n \left(1+\frac{s}{\gamma_i(P)}\right), 
\end{equation*}
where $-\gamma_i(P)\in \C$, $1\leq i\leq n$,  are the roots of the
Ehrhart polynomial $\LE(s,P)$. In particular, for their geometric and
arithmetic mean we have 
\begin{equation}
  \left(\prod_{i=1}^n \gamma_i(P)\right)^{1/n} = \left(1/\vol(P)\right)^{1/n},\quad 
  \frac{1}{n}\sum_{i=0}^n \gamma_i(P) = \frac{1}{n} \frac{\LE_{n-1}(P)}{\vol(P)}.
\label{eq:coeffroots}
\end{equation}

Here we are interested in geometric interpretations of the roots and
in their size. Since the volume of lattice polytopes
without interior lattice points might be arbitrary large for $n\geq 3$
the norm of the roots $|\gamma_i(P)|$  might be arbitrary small
(cf.~\eqref{eq:coeffroots}). On the other hand, we know that the
volume of a lattice polytope with $l\geq 1$ interior  lattice points 
is bounded 
(cf.~\cite{hensley:lattice_volume,  lagariasziegler:volume,
pikhurko:lattice_volume, zpw:82}) by a constant depending only on $l$ and
$n$. Thus, up to unimodular transformations, there are only finitely
many different of those lattice polytopes and in this case  
$|\gamma_i(P)|$ can not be too small. 
Therefore, it seems to be reasonable to distinguish lattice polytopes
with or without 
interior lattice points, and we define 
 \begin{definition} 
For $l\in\N$ let $\mathcal{P}^n(l)$ be the set of lattice
  polytopes $P\in \mathcal{P}^n$ having exactly $l$ interior lattice points, i.e., $\LE(\inter(P))=l$.
Moreover, 
the set of all $0$-symmetric lattice polytopes $P\in\mathcal{P}^n$ 
is denoted by $\mathcal{P}_o^n$.
\end{definition}
As already mentioned above, it was shown by Pikhurko \cite{pikhurko:lattice_volume} that for
$P\in\mathcal{P}^n(l)$, $l\geq 1$, 
\begin{equation*}
   \vol(P) \leq  \mathrm{c}_n\,l, 
\end{equation*}  
for a constant $\mathrm{c}_n$ depending only on $n$. Hence we get  
\begin{equation*}
 \left(\prod_{i=1}^n \gamma_i(P)\right)^{1/n} \geq (\mathrm{c}_n)^{-1/n} l^{-1/n}.
\end{equation*}
Candidates of lattice polytopes $P\in\mathcal{P}^n(l)$, $l\geq 1$, of
maximal volume are certain simplices $T(n,l)\in \mathcal{P}^n(l)$, introduced by Perles, Wills and Zaks
\cite{zpw:82} with $\vol(T(n,l))\geq (l+1)/n!\, 2^{2^{n-1}}$.

In order to present a lower bound on the volume  in terms of the
number of interior lattice points we define for
$l\in\N$ the simplices 
\begin{equation*}
  S_n(l)=\conv\Big\{e_1,\dots,e_n,-l\,\sum_{i=1}^n e_i\Big\},
\end{equation*} 
where $e_i$ denotes the $i$-th unit vector. Observe that $\LE(\inter(S_n(l)))=l$ and $\vol(S_n(l))=(n\,l+1)/n!$. 
\begin{theorem} Let $P\in\mathcal{P}^n$. Then 
\begin{equation}
   \vol(P) \geq \frac{n\,\LE(\inter \,P)+1}{n!},
\label{eq:non-symmetric}
\end{equation}
and the bound is best possible for any number of interior lattice points. For $\LE(\inter P)=1$ equality holds if and only if $P$ is unimodular isomorphic  to
the simplex $S_n(1)$. 
\label{thm:non-symmetric}
\end{theorem} 
The theorem above implies that for $P\in \mathcal{P}^n(l)$ the geometric mean of the roots is bounded from above by 
\begin{equation*}
 \left(\prod_{i=1}^n \gamma_i(P)\right)^{1/n} \leq (n!)^{1/n} (n\,l+1)^{-1/n}.
\end{equation*} 
In the $2$-dimensional case Theorem \ref{thm:non-symmetric} is a
direct consequence of Pick's identity $G(P)=\vol(P)+\frac{1}{2}\LE(\bd
P)+1$, where $\bd P$ denotes the
boundary  of $P$ \cite{pick:lattice}. In particular, equality is
attained in \eqref{eq:non-symmetric} iff
$P$ is a lattice triangle whose vertices are the only lattice points
contained in the boundary. This also shows that for 
$\LE(\inter P)>1$ the
extremal cases  in \eqref{eq:non-symmetric} are not necessarily unimodular
equivalent. We remark, however, 
 that all extremal cases have the same Ehrhart
polynomial. 
\begin{proposition} Let $P\in\mathcal{P}^n(l)$, $l\geq 1$, with
  $\vol(P)=(n\,l+1)/n!$. Then $\sa_i(P)=\sa_i(S_n(l))=l$, $1\leq i\leq
  n$.
\label{prop:non-symmetric} 
\end{proposition}

For $0$-symmetric lattice polytopes $P\in\mathcal{P}_o^n$ there is a classical upper bound on the
volume due to Blichfeldt  and van der Corput
(cf.~\cite[p.~51]{gruber-lekkerkerker:87})
\begin{equation*}
 \vol(P)\leq 2^{n-1}\, \left(\LE(\inter P)+1\right). 
\end{equation*} 
Lattice boxes 
\begin{equation*}
 Q_n(2\,l-1)=\{x\in\R^n : |x_1|\leq l,\,|x_i|\leq 1,\, 2\leq i\leq n\}
 \in\mathcal{P}^n_o,\quad l \in\N ,
\end{equation*}
show that the bound is tight. As an analogue to Theorem
\ref{thm:non-symmetric} in the $0$-symmetric case we conjecture 
\begin{conjecture} Let $P\in\mathcal{P}_o^n$. Then 
\begin{equation*}
  \vol(P)\geq \frac{2^{n-1}}{n!}\, \left(\LE(\inter P)+1\right).
\end{equation*}
\label{conj:symmetric}
\end{conjecture}
Again for $n=2$ the inequality follows immediately from  Pick's identity
and  the inequality is tight for
any parallelogram whose vertices are the only lattice points on the
boundary. It seems to be quite likely that certain crosspolytopes, i.e.,
$P\in\mathcal{P}_o^n$ with $2\,n$ vertices, are the extremal cases for
the inequality above; for the family of $0$-symmetric crosspolytopes we can prove the conjecture.
\begin{proposition} For $P\in\mathcal{P}_o^n$ with $2\,n$
  vertices Conjecture \ref{conj:symmetric} is true.
\label{prop:crosspolytopes}
\end{proposition}   
One way to prove  that proposition is based on the following lemma which
might be of some interest in its own.
\begin{lemma} Let $P\in\mathcal{P}_o^n$ with $2\,n$ vertices. Then 
\begin{equation*}
   a_i(P)+a_{n-i}(P) \geq \binom{n}{i}\left(a_n(P)+a_0(P)\right),\quad
   i=0,\dots,n.
\end{equation*} 
\label{lem:favourite}  
\end{lemma} 
 Observe that on account of \eqref{eq:stanley_coefficients} Lemma
 \ref{lem:favourite} implies Proposition \ref{prop:crosspolytopes}. 
As a side effect of the proof of that lemma we get 
\begin{remark}  Let $P\in \mathcal{P}_o^n$. Then 
  $\sa_i(P)\geq \binom{n}{i}$ for $\,0\leq i\leq n$.
\label{rem:symmetric}
\end{remark}
The lower bounds on $\sa_i(P)$  in the remark above also follow from a much deeper and much more general
result of Stanley \cite{stanley:symmetric} on the $h$-vector of ``symmetric'' Cohen-Macaulay simplical
complexes in conjunction with a result of Betke and McMullen \cite{betke-mcMullen:85} relating the
coefficients $\sa_i(P)$ with the $h$-vector of a 
triangulation of the polytope. Here we give a quite elementary proof
which follows the method
presented by Beck and Sottile in \cite{beck_sottile:irrational}.

The regular
 unit crosspolytope $C_n^\star=\conv\{\pm e_i : 1\leq i\leq n\}$ plays
 a special role in the context of the roots of Ehrhart polynomials. To
 our knowledge it was firstly shown by 
Kirschenhofer, Pethoe, and Tichy
\cite{kirschenhofer_pethoe_tichy:counting_polynomials} that the real part of
$-\gamma_i(C_n^\star)$ is equal to $-1/2$, $1\leq i\leq n$. This was independently
proven by Bump et al. in \cite[Theorem 4]{bump:crosspolytope} and follows also from
a more general result of Rodriguez-Villegas \cite{rodriguez:zeros}. In
\cite[Open problem 2.42]{beckrobins:discrete_continously} the authors
ask for other classes of lattice polytopes such that all roots of their Ehrhart  polynomials have real part $-1/2$. Since $C_n^\star$  has minimal
volume among all $0$-symmetric lattice polytopes an obvious candidate
is the simplex $S_n(1)$ in the non-symmetric case (cf.~Theorem
\ref{thm:non-symmetric}).  

\begin{theorem} All roots of the polynomial $\LE(s,S_n(1))$ have real
  part $-1/2$. If $\alpha_n$ is a root of $\LE(s,S_n(1))$ with
  maximal norm, then \begin{equation*}
    \left|\alpha_n+\frac{1}{2}\right|=\frac{n(n+2)}{2\,\pi}+o(n), 
\end{equation*} 
as $n$ tends to infinity.
\label{thm:roots}
\end{theorem}

In a recent paper Braun \cite{braun:ehrhart_roots} proved that the 
roots of an Ehrhart polynomial lie inside the disc with center $-1/2$
and radius $n(n-1)/2$. The above theorem shows that this bound is
essentially  tight and improves on the former best known bound of
order $n$ \cite[Theorem 1.3]{becketall:roots_ehrhart}.

It seems to be quite likely that $\LE(s,S_n(1))$ possesses  the roots
of  maximal norm among all Ehrhart polynomials of polytopes with  interior points. In the case $n=2$ this follows from \cite[Theorem 2.2]{becketall:roots_ehrhart} and for a verification of this statement in the 3-dimensional case see Theorem \ref{thm:3polytopes}.

Looking at geometric properties of lattice polytopes $P$ whose roots
have all real part $-1/2$ leads immediately to the class of reflexive
lattice polytopes. Here $P\in\mathcal{P}^n$ with $0\in\inter P$ is called reflexive if 
\begin{equation*}
   P^\star=\{y\in\R^n : x\,y \leq 1,\text{ for all } x\in P\}\in\mathcal{P}^n, 
\end{equation*}
i.e., the polar polytope is again a lattice polytope. They play an
important role in toric geometry since they are in one-to-one
correspondence with Gorenstein toric Fano varieties. Reflexive polytopes have been extensively studied and exhibit many surprising properties (cf.~\cite{batyrev:reflexive, hibi:92, nill:thesis} and the references within). 
 In particular, Hibi \cite{hibi:92} showed that the coefficients
 $\sa_i(P)$ of a lattice polytope are symmetric, i.e.,
 $\sa_i(P)=\sa_{n-i}(P)$,  if and only if $P$ is reflexive. 
  Kreuzer and Skarke \cite{kreuzer_skarke:reflexive3,
   kreuzer_skarke:reflexive4} classified all reflexive polytopes in
 dimensions $\leq 4$. For $n=2,3,4$ there are respectively $16$; $4,319$ and $473,800,776$ reflexive polytopes (up to unimodular equivalence).

\begin{proposition} Let $P\in\mathcal{P}^n$. If all roots of
  $\LE(s,P)$ have real part $-1/2$ then, up to an unimodular translation, 
  $P$ is a reflexive polytope
  of volume  $\leq 2^n$.
\label{prop:roots_reflexive}
\end{proposition} 

It is easy to check that for $n\leq 3$ the converse is also true but
not for $n\geq 4$. All in all, for $n\leq 4$ we the following
characterization  
\begin{proposition} Let $P\in\mathcal{P}^n$ be a reflexive polytope. Then all roots of
  $\LE(s,P)$ have real part $-1/2$ 
\hfill \renewcommand{\labelenumi}{\rm \roman{enumi})} 
\begin{enumerate}
\item iff $\vol(P)\leq 2^n$ and $n\leq 3$,  
\item iff  $(\LE(P)-1-4\,\vol(P))^2\geq 16\,\vol(P)$, $\LE(P)\leq 9\,\vol(P)+18$ and $n=4$. 
\end{enumerate} 
\label{prop:ref_roots}
\end{proposition}   
A classification of the roots of $2$-dimensional lattice polygons is
    given in the papers \cite[Theorem 2.2]{becketall:roots_ehrhart}
    and \cite[Theorem 1.9]{henkschuermannwills:ehrhart}.  
For $n=3$ we know less  and the basic properties are subsumed in the next theorem. For more
    detailed properties of Ehrhart polynomials of $3$-dimensional lattice polytopes we refer to section \ref{sec:3polytopes}.
\begin{theorem} The roots of the Ehrhart polynomials of
  $3$-dimensional lattice polytopes are contained in 
\begin{equation*}
 [-3,-1]\cup\left\{a+\mathrm{i}\,b : -1\leq a <1,\, a^2+b^2\leq 3\right\}
\end{equation*}
and the bounds on $a$ and $a^2+b^2$  are tight. For $P\in\mathcal{P}_3(l)$, $l\geq 1$, the upper
bound $\sqrt{3}$  on the norm of the complex roots  is  
  only attained by the roots of the Ehrhart polynomial of the simplex $S_3(1)$.
\label{thm:3polytopes}
\end{theorem}

The paper is organized as follows. In Section 2 we prove Theorem
\ref{thm:non-symmetric}, Theorem \ref{thm:roots} and what we know in the $0$-symmetric case
regarding Conjecture \ref{conj:symmetric}. Section 3 deals with
reflexive polytopes and Ehrhart polynomials whose roots have all 
real part $-1/2$. Section 3 studies the Ehrhart polynomials and their
roots for $3$-dimensional lattice polytopes.

\section{Volume and interior lattice points}

The proof of Theorem \ref{thm:non-symmetric} is based on a subdivision
of $P$ with respect to the interior lattice points contained in $P$.
\begin{proof}[Proof of Theorem \ref{thm:non-symmetric}] Let $l=\LE(\inter(P))$. If $l=0$ there is nothing to show since any lattice polytope has at least volume $1/n!$. So let $l>0$ and let $y_1,\dots,y_l$ be the interior lattice points of $P$. Obviously, it suffices to show that $P$  can be subdivided with the points $y_1,\dots,y_k$ into at least $n\,k+1$ lattice polytopes for $k=1,\dots,l$.  

First we build the convex hulls of $y_1$ with all facets of $P$ yielding at least $n+1$ lattice polytopes.  So let us assume that we have already dissected $P$ into $P_1,\dots,P_{n\,k+1}$ lattice polytopes and let $y_{k+1}$ be contained in the relative interior of a $j$-dimensional face  of $P_1$, say. Since $y_{k+1}$ is an interior point it is also contained in the relative interior  of a $j$-face  of at least $n-j$ further polytopes $P_2,\dots,P_{r}$, say, $r\geq n-j+1$. Subdividing each $P_s$ by building the convex hull of $y_{k+1}$ with all facets of $P_s$ not containing $y_{k+1}$ gives at least $j+1$ new polytopes for each $P_s$, $s=1,\dots,r$. Thus this new subdivision of $P$  consists of at least 
\begin{equation*}
   r\cdot(j+1)+ n\,k+1-r \geq (n-j+1)j+n\,k+1 
\end{equation*}
lattice polytopes. Since $j\geq 1$ this number is at least  $n\,(k+1)+1$. 

The simplices $S_n(l)$ show that the bound is attained for any number of interior lattice points and the proof above shows that equality can only be achieved by simplices. So let us assume that we have a lattice simplex $S$ with only one interior lattice point $y_1$ and equality in equation \eqref{eq:non-symmetric}. Without loss of generality let  $y_1=0$ and let $v_1,\dots,v_{n+1}$ be the vertices of $S$. Let $F_i$ be the facet of $S$ not containing $v_i$, $1\leq i\leq n+1$. Subdividing $S$ into the $n+1$ simplices $\conv\{0,F_i\}$, $1\leq i\leq n+1$, gives 
\begin{equation*}
    \frac{n+1}{n!}=\vol(S)=\sum_{i=1}^{n+1} \vol(\conv\{0,F_i\}).
\end{equation*} 
Since $\vol(\conv\{0,F_i\})\geq 1/n!$ we must have $\vol(\conv\{0,F_i\})=1/n!$, or equivalently, any choice of $n$ vectors out of the vertices form a basis of the lattice $\Z^n$. Thus, up to an unimodular linear transformation, we may assume $v_i=e_i$, $1\leq i\leq n$, and the absolute value of each coordinate of $v_{n+1}$ is $1$. Finally, since $0$ is contained in the interior of $S$ we must have $v_{n+1}=(-1,\dots,-1)^\intercal$.  
\end{proof} 

\smallskip
We remark that inequality \eqref{eq:non-symmetric} can also be deduced
from a result of Hibi \cite{Hibi:lower_bound} where it is shown that 
\begin{equation}
  \sa_i(P)\geq \sa_1(P), \, 1\leq i\leq n-1,
\label{eq:hibi} 
\end{equation} 
for $P\in\mathcal{P}^n(l)$ with $l\geq 1$. Together with
\eqref{eq:stanley_coefficients} 
this implies \eqref{eq:non-symmetric}. 

From \eqref{eq:hibi} we also get  
Proposition \ref{prop:non-symmetric}, i.e.,  
the uniqueness of the Ehrhart polynomials of $P\in\mathcal{P}^n(l)$,
$l\geq 1$, with minimal volume. By \eqref{eq:hibi} 
we know $\sa_i(P)\geq l$ for $1\leq i\leq n$
(cf.~\eqref{eq:stanley_coefficients}) and since $P$ has minimal volume
we also have 
\begin{equation*}
n!\vol(S_n(l))= 1+n\,l = n!\vol(P)=\sa_0(P)+\sa_1(P)+\ldots+\sa_n(P).
\end{equation*}
Hence $\sa_i(P)=\sa_i(S_n(l))=l$, $1\leq i\leq n$,  
and the Ehrhart polynomial of $P\in\mathcal{P}_n(l)$ with minimal
volume is uniquely determined. 
\smallskip


We believe that the crosspolytopes 
\begin{equation*}
  C_n^\star(2\,l-1)=\conv\{\pm l\,e_1,\pm\,e_2,\dots,\pm\,e_n\},\quad l\geq 1,
\end{equation*} 
with $2l-1$ interior lattice points form the $0$-symmetric counterpart to the simplices $S_n(l)$, i.e., they have minimal volume among all $0$-symmetric polytopes with $2l-1$ interior lattice points. In \cite[Theorem 2.6]{beckrobins:discrete_continously} it is shown that the coefficients $\sa_i(P)$ of a bipyramid $P=\conv\{Q,\pm e_n\}$, where $Q$ is an $(n-1)$-dimensional lattice polytopes embedded in the hyperplane $\{x\in\R^{n-1} :x_n=0\}$ satisfy the recursion $\sa_i(P)=\sa_i(Q)+\sa_{i-1}(Q)$. Hence we conclude  
\begin{equation*}
   \sa_i(C_n^\star(2\,l-1))=\binom{n}{i}+\binom{n-1}{i-1}(2\,l-2), 
\end{equation*} 
and so 
\begin{equation*}
   \sa_i(C_n^\star(2\,l-1))+\sa_{n-i}(C_n^\star(2\,l-1))
   =\binom{n}{i}\left(\sa_n(C_n^\star(2\,l-1))+1\right), \quad 0\leq i\leq n.
\end{equation*} 
Lemma \ref{lem:favourite} shows that $\binom{n}{i}\left(\sa_n(P)+1\right)$ is a lower bound on $\sa_i(P)+\sa_{n-i}(P)$ for any lattice crosspolytope $P$. 

\begin{proof}[Proof of Lemma \ref{lem:favourite}] Let
  $P=\conv\{\pm\,v_j : 1\leq j\leq n\}$ be a lattice crosspolytope in
  $\R^n$.  For any of the $2^n$ subsets 
  $W_k=\{w_1,\dots,w_n\}$ with $w_j\in\{v_j,-v_j\}$  we consider 
  the simplicial cone 
  $C_k=\cone\{(w_1,1)^\intercal,\dots,(w_n,1)^\intercal,e_{n+1}\}\subset\R^{n+1}$ 
  and the open parallelepiped 
\begin{equation*}
  Q_k =\Big\{\sum_{j=1}^n\lambda_j \binom{w_j}{1} + \lambda_{n+1}\,e_{n+1}: 0<\lambda_j<1,\,1\leq j\leq n+1\Big\}. 
\end{equation*}
The cones $C_k$ form a triangulation of the cone   
$C=\cone\{(\pm\,v_j,1)^\intercal : 1\leq j\leq n\}$.
In a recent paper Beck and Sottile\cite{beck_sottile:irrational}
introduced a new method for ''calculating'' the numbers
$\sa_i(\cdot)$ of an arbitrary lattice polytope. In order to apply 
their approach we choose a vector
$s=(s_1,\dots,s_{n+1})^\intercal\in\R^{n+1}$ such that  
$C\cap\Z^{n+1}= (s+C)\cap\Z^{n+1}$ 
and none of the shifted cones $s+C_k$ contains a lattice point on its
boundary. Obviously, we must have $s_{n+1}<0$ and the vector $s$ can be chosen
arbitrarily short. For $i=0,\dots,n$ we denote by  
\begin{equation*}
  \alpha_i(s+Q_k)= \#\left\{ (s+Q_k)\cap\{z\in\Z^{n+1} : z_{n+1}=i\}\right\}
\end{equation*} 
the number of lattice points in $s+Q_k$ having last coordinate
$i$. Then  we have (see \cite[Proof of Theorem 2]{beck_sottile:irrational}, \cite[Proof of Theorem 3.12]{beckrobins:discrete_continously})
\begin{equation}
   \sa_i(P) = \sum_{k} \alpha_i(s+Q_k).
\label{eq:summe}
\end{equation} 
In particular we have $\sa_n(P)$ many lattice points with last
coordinate $n$ contained in the parallelepipeds $s+Q_k$.  Let
$(w,n)^\intercal$ be one of them and let it be given by 
\begin{equation}
 \binom{w}{n} = s + \sum_{j=1}^n \lambda_j\,\binom{w_j}{1} +\lambda_{n+1}\, e_{n+1}
\label{eq:one}
\end{equation} 
Now we fix an $i\in\{1,\dots,n-1\}$. For a subset
$I\subset\{1,\dots,n\}$ of cardinality $i$ we denote by $I^c$ its
complement and let 
\begin{equation*}
  \lambda_I:= \lambda_{n+1}+2\sum_{j\in I} (\lambda_j-1).
\end{equation*} 
With this notation we may write 
\begin{equation*}
\begin{split} 
f(w,I) & :=\binom{w}{n}-\sum_{j\in I} \binom{w_j}{1} \\
& =  s +\sum_{j\in I}
(1-\lambda_j)\binom{-w_j}{1}+\sum_{j\in I^c}\lambda_j\,\binom{w_j}{1}
+ \lambda_I\, e_{n+1} 
\end{split}
\end{equation*}
Hence, if the scalar $\lambda_I$ is positive then the lattice point $f(w,I)$ 
is contained in some $s+Q_{k^\prime}$, say, and therefore, it 
contributes to $\sa_{n-i}(P)$ (cf.~\eqref{eq:summe}).

Since $s_{n+1}<0$ we have 
$\sum_{j=1}^{n+1} \lambda_j > n$ (cf.~\eqref{eq:one}) and so we get
either $\lambda_I>0$ or $\lambda_{I^c}>0$.
 In other words, either $f(w,I)$ contributes to $\sa_{n-i}(P)$ or
 $f(w,I^c)$ contributes to $\sa_{i}(P)$. Since this argument works for
 any subset $I$ of cardinality $i$  the lattice point $(w,n)^\intercal$ ''produces'' in this way a contribution of $\binom{n}{i}$ to the sum $\sa_i(P)+\sa_{n-i}(P)$. 

Next we have to check that for two different points 
\begin{equation*} 
\binom{w}{n} = s + \sum_{j=1}^n \lambda_j\,\binom{w_j}{1} +\lambda_{n+1}\, e_{n+1} \text{ and }  
\binom{\widetilde{w}}{n} = s + \sum_{j=1}^n \mu_j\,\binom{\widetilde{w}_j}{1} +\mu_{n+1}\, e_{n+1}, 
\end{equation*}
the lattice points 
$f(w,I)$ and $f(\widetilde{w},\widetilde{I})$, $\#I=\#\widetilde{I}=i$,  are also different, provided both
of them  contribute to $\sa_{n-i}(P)$. 
Suppose the opposite, i.e., $f(w,I)=f(\widetilde{w},\widetilde{I})$. Since both of
them contribute to $\sa_{n-i}(P)$ the two points $f(w,I)$,
$f(\widetilde{w},\widetilde{I})$ lie in the same cone $s+C_{k^\prime}$, say, and
since any lattice point in $s+C$ is contained in exactly one of the simplicial cones $s+C_k$ we conclude 
\begin{equation*}
\{-w_j: j\in I\}\cup\{w_j:j\in I^c\} = 
\{-\widetilde{w}_j: j\in \widetilde{I}\}\cup\{\widetilde{w}_j : j\in {\widetilde{I}}^c\}.
\end{equation*}
If $\{-w_j: j\in I\}=\{-\widetilde{w}_j: j\in \widetilde{I}\}$ then we must also have $\{w_j:j\in I^c\}=\{\widetilde{w}_j : j\in {\widetilde{I}}^c\}$.
Each point in a simplicial cone, however, has an unique representation
with respect to the generators and so we get the contradiction
$(w,n)^\intercal=(\widetilde{w},n)^\intercal$. 
Therefore, we may assume that there exists a $j_1\in I\cap {\widetilde{I}}^c$ and  a $j_2\in \widetilde{I}\cap I^c$.
 Thus 
$1-\lambda_{j_1}=\mu_{j_1}$ and $1-\mu_{j_2}=\lambda_{j_2}$ and so 
\begin{equation} 
    \mu_{j_1} +\mu_{j_2}+\lambda_{j_1}+\lambda_{j_2}=2.
\label{eq:two}
\end{equation} 
On the other hand, since $\sum_{i=1}^{n+1}\lambda_i$,
$\sum_{i=1}^{n+1}\mu_i  > n$ and $\lambda_i$, $\mu_i<1$ we have  $\lambda_{j_1}+\lambda_{j_2},\mu_{j_1} +\mu_{j_2}>1$ contradicting \eqref{eq:two}.

So far we have shown that 
\begin{equation} 
 \sa_i(P)+\sa_{n-i}(P)\geq \binom{n}{i}\sa_n(P).
\label{eq:sofar}
\end{equation}
 Now there is one special point $\binom{\widehat{w}}{n}$
which contributes to $\sa_{i}(P)$ as well as  to $\sa_{n-i}(P)$. 
Since the origin $0\in\R^{n+1}$ is contained in
one of the cones $s+C_k$, say, we can find a representation of the form 
\begin{equation*}
  0 = s+ \sum_{j=1}^n \mu_j \binom{w_j}{1} + \mu_{n+1}\,e_{n+1},
\end{equation*} 
$\mu_i>0$. Choosing the vector $s$ sufficiently small we may assume that 
\begin{equation}
   \mu_{n+1} + 2\,\sum_{j=1}^n \mu_j < 1.
\label{eq:three} 
\end{equation} 
Hence the vector
\begin{equation*}  
\binom{\widehat{w}}{n}=\sum_{j=1}^n \binom{-w_j}{1} = s + \sum_{j=1}^n (1-\mu_j)\binom{-w_j}{1}+\left(\mu_{n+1} + 2\,\sum_{j=1}^n \mu_j \right)e_{n+1} 
\label{eq:vectorhat}
\end{equation*}
is contained in some $s+Q_{k^\prime}$, say. On account of \eqref{eq:three} the vectors $f(\widehat{w},I)$ and $f(\widehat{w},I^c)$ are contained in
some of these parallelepipeds  for all subsets $I\subset\{1,\dots,n\}$ of
cardinality $i$.  Thus the vector $(\widehat{w},n)^\intercal$
gives a contribution of $\binom{n}{i}$ to $\sa_i(P)$ and
to $\sa_{n-i}(P)$.  Together with \eqref{eq:sofar} this proves the lemma.
\end{proof}

For the proof of the inequalities in Remark \eqref{rem:symmetric}
we just observe that the last part of the proof above where the vector
$(\widehat{w},n)^\intercal$ is considered, in particular
implies that $\sa_i(P)\geq \binom{n}{i}$ for any lattice
crosspolytope. Now any $n$-dimensional $0$-symmetric lattice polytope
$\tilde{P}$ 
contains a $0$-symmetric lattice crosspolytopes $P$ and by Stanley's
Monotonicity Theorem \cite{stanley:monoton} (see also
\cite{beck_sottile:irrational}) we have $\sa_i(\tilde{P})\geq\sa_i(P)$.   

Next we come to the roots of the polynomial $\LE(s,S_n(1))$ which on
account of Proposition \ref{prop:non-symmetric} is given by
\begin{equation}
 \LE(s,S_n(1)) =\sum_{i=0}^n \binom{s+n-i}{n} = \binom{s+n+1}{n+1}-\binom{s}{n+1}.
\label{eq:simplex} 
\end{equation} 
 
\begin{proof}[Proof of Theorem \ref{thm:roots}] One way to see that
  all roots of $\LE(s,S_n(1))$ have real part $-1/2$ is to apply a
  theorem of Rodriguez-Villegas \cite{rodriguez:zeros}. In our setting it says that  if all roots of the polynomial $f(s,P)=\sum_{i=0}^n \sa_i(P)\,s^i$
  lie on the unit circle then all roots of $\LE(s,P)$ have real part
  $-1/2$. In our case we have $f(s,S_n(1))=\sum_{i=0}^n s^i$ and
  so the norm of each root of that polynomial is $1$.

  Now let
  $s_0=-1/2+\mathrm{i}\,b=r_0\,\mathrm{e}^{\mathrm{i}\,\alpha_0}$ be a
  point on the line with real part $-1/2$ where 
  we assume  $b\geq 0$. Furthermore, for  
  $m=1,\dots,n$ let
  $s_0-m=r_m\,\mathrm{e}^{\mathrm{i}\,\alpha_m}$. Since
  $|s_0-m|=|s_0+m+1|$, $m=0,\dots,n$, we also know that
  $s_0+m+1=r_m\,\mathrm{e}^{\mathrm{i}\,(\pi-\alpha_m)}$. From
  the right hand side of \eqref{eq:simplex} we conclude that $s_0$ is
  a root of $\LE(s,S_n(1))$ if and only if  
\begin{equation*}
  (s_0+n+1)\,(s_0+n)\cdot\ldots\cdot(s_0+1)=s_0\,(s_0-1)\cdot\ldots\cdot (s_0-n).
\end{equation*} 
Substituting the polar representations leads to 
\begin{equation*}
   (-1)^{n+1} = \mathrm{e}^{\mathrm{i}\,(2\,\alpha_0+2\,\alpha_1+\ldots + 2\,\alpha_n)}
\end{equation*}
Replacing the angle $\alpha_m$ by $\pi/2+\overline{\alpha_m}$, $\overline{\alpha_m}\in (0,\pi/2]$,  gives 
$
   1 = \mathrm{e}^{\mathrm{i}\,(2\,\overline{\alpha_0}+2\,\overline{\alpha_1}+\ldots + 2\,\overline{\alpha_n})}
$ and thus we must have 
\begin{equation*}
   \overline{\alpha_0}+\overline{\alpha_1}+\ldots +
   \overline{\alpha_n} = k\,\pi.
\end{equation*} 
for an integer $k\in\{1,\dots,\lfloor (n+1)/2\rfloor\}$. Observe, that
we have assumed $b\geq 0$. By
construction we have $\cot \overline{\alpha_m}  = 2\,b/(2\,m+1)$, $m=0,\dots,n$, 
and so we get that $s_0=-1/2+\mathrm{i}\,b$ is a root of
$\LE(s,S_n(1))$ if  and only if 
\begin{equation*}
 h(b):= \sum_{m=0}^n \cot^{-1}\left(\frac{2\,b}{2\,m+1}\right)
 \in\{\pi,2\,\pi,\dots, \lfloor(n+1)/2\rfloor\,\pi\},
\end{equation*} 
where we require $\cot^{-1}()\in(0,\pi/2]$. 
Since $h(b)$ is a monotonously decreasing function in $b$ the
imaginary part $b_n$ 
of the  root of maximal norm is determined by the equation $h(b_n)=\pi$. 
Since $\cot^{-1}(t)=\tan^{-1}(1/t)$ ``the inverse'' of the cotangent has the power series representation
$\cot^{-1}(t)=\sum_{k=0}^\infty (-1)^k/(2k+1)\,(1/t)^{2k+1}$ for $t>1$. So
we have $1/t>\cot^{-1}(t)>1/t-1/(3t^3)$. Hence for $b>n+1/2$
we may write 
\begin{equation*}
              \frac{1}{b}\frac{(n+1)^2}{2} > h(b) >
              \frac{1}{b}\frac{(n+1)^2}{2} - \mathrm{c}\,\frac{n^4}{b^3}
\end{equation*}  
for a suitable constant $\mathrm{c}$. Thus $b_n= n(n+2)/(2\pi)+o(n).$
\end{proof}

\section{Reflexive polytopes }
As mentioned in the introduction reflexive polytopes are of
particular interests in many different branches of mathematics and
have a lot of nice geometric properties.  Some of them 
 are collected in the
following lemma for which we refer to \cite{batyrev:reflexive, hibi:92}.
\begin{lemma} Let $P\in\mathcal{P}^n$ with $0\in\inter(P)$. Then $P$
  is relexive if and only if \renewcommand{\labelenumi}{\rm \roman{enumi})} 
\begin{enumerate}
\item $P^\star\in \mathcal{P}^n$.
\item $\sa_i(P)=\sa_{n-i}(P)$, $0\leq i\leq n$.
\item $\LE(k\,P)=\LE((k+1)\,\inter(P))$ for $k\in\N$.
\item $\vol(P)=(n/2)\,\LE_{n-1}(P)$, i.e., the origin lies in an
  adjacent lattice hyperplane to any facet.
\end{enumerate} 
\label{lem:relexive_properties}
\end{lemma}
In particular, the origin is the only interior lattice point of a
reflexive polytope, and reflexive polytopes are precisely those lattice polytopes satisfying the functional equation: 
\begin{equation*}
\LE(s,P) =(-1)^n\LE(-(1+s),P),\, s\in\C.
\end{equation*}
Hence in any odd dimension the Ehrhart polynomials of reflexive polytopes have the real root $-1/2$.

Now let $P\in\mathcal{P}^n$ be a lattice polytope such that the real
part of all roots $-\gamma_i(P)$ of its Ehrhart polynomial is -$1/2$. Then from \eqref{eq:coeffroots} we immediately  get
\begin{equation*} 
  \frac{n}{2}\LE_{n-1}(P)=\vol(P)\leq 2^n 
\end{equation*}  
which by Lemma \ref{lem:relexive_properties} iv) verifies Proposition \ref{prop:roots_reflexive}.

In dimension 2 any lattice polygon $P$ whose only interior lattice point
is the origin is reflexive and its Ehrhart polynomial is given by (cf.~\eqref{eq:ehrhart_series}, \eqref{eq:stanley_coefficients})
\begin{equation*}
\LE(s,P)=\vol(P)\,\left( s^2+s+\frac{1}{\vol(P)}\right).
\end{equation*}
Thus all roots have real part $-1/2$ if and only if
$\vol(P)\leq 4$. Among the well known 16 reflexive
polytopes in $\R^2$ (cf.~e.g.~\cite{Poonen_Rodriguez-Villegas:reflexive}) 
there is only one with volume bigger than $4$, 
namely the simplex $S=-(1,1)^\intercal+\conv\{0,3\,e_1,3\,e_2\}$ of
volume $9/2$. By Theorem \ref{thm:non-symmetric} we know that 
the reflexive polygon of minimal volume is $S_2(1)$ of volume $3/2$. 
Hence the Cartesian product $S\times S_2(1)$ is an example of a
$4$-dimensional reflexive polytope of volume less than $2^n$ ($n=4$), but not
all roots of its Ehrhart polynomial have real part $-1/2$.   

\begin{proof}[Proof of Proposition \ref{prop:ref_roots}] First we check that all  roots of the Ehrhart polynomial  of a 3-dimensional reflexive polytope $P$ have real part $-1/2$ if and only if its volume is not bigger than $8$.  By Lemma \ref{lem:relexive_properties} ii) we have $\sa_1(P)=\sa_2(P)$ and so (cf.~\eqref{eq:ehrhart_series}, \eqref{eq:stanley_coefficients})
\begin{equation*}
\begin{split}
 \LE(s,P) & = \frac{1}{6}\left[(2\sa_1(P)+2)\,s^3+(3\sa_1(P)+3)\,s^2+(13\sa_1(P)+1)\,s+6\right] \\
&=\vol(P)\left[s^3+\frac{3}{2}s^2+\left(\frac{1}{2}+\frac{2}{\vol(P)}\right)\,s+\frac{1}{\vol(P)}\right] \\
&= \vol(P)\left[ \left(s+\frac{1}{2}\right)\left(s^2+s+\frac{2}{\vol(P)}\right)\right]. 
\end{split}
\end{equation*}   
Hence all roots have real part $-1/2$ iff $\vol(P)\leq 8$. 

Now let $P$ be a $4$-dimensional reflexive polytope. Again by Lemma  \ref{lem:relexive_properties} we have $\sa_1(P)=\sa_3(P)$ and so we find 
\begin{equation*}
\begin{split}
\LE(s,P)  & = \frac{1}{24}\Big[ (2\sa_1(P)+\sa_2(P)+2)\,s^4+(4\sa_1(P)+2\sa_2(P)+4)\,s^3 \\ &\quad\quad\quad+(10\sa_1(P)-\sa_2(P)+46)\,s^2+(8\sa_1(P)-2\sa_2(P)+44)\,s+24\Big] \\
& = \vol(P)\Big[ s^3+2\,s^3 + \left(2\,\mu+1\right)\,s^2 +  \left(2\mu\right)\,s+\frac{1}{\vol(P)}\Big],
\end{split}
\end{equation*}
where $\mu=(1+(1/4)\sa_1(P))/\vol(P)-1$. Further we set $\beta=1/\vol(P)$ and obtain 
\begin{equation*}
 \LE(s,P) = \vol(P)\Big[ \left(s^2+s+\mu+\sqrt{\mu^2-\beta}\right)\cdot\left(s^2+s+\mu-\sqrt{\mu^2-\beta}\right)\Big].
\end{equation*}
Thus all roots have real part $-1/2$ if and only if $\mu^2\geq \beta$ and $\mu-\sqrt{\mu^2-\beta}\geq 1/4$.  The first condition translates into $(2+(1/2)\sa_1(P)-2\,\vol(P))^2\geq 4\,\vol(P)$ and the seond becomes $2\,\sa_1(P)\leq 9\,\vol(P)+8$. Since $\sa_1(P)=\LE(P)-5$ we get the inequalities stated in Proposition \ref{prop:ref_roots}. 
\end{proof}

Thanks to the classification  of Kreuzer and Skarke (cf.~\url{http://hep.itp.tuwien.ac.at/~kreuzer/CY/}) one can check that among the\
 4319 reflexive polytopes in dimension 3 only 64 have volume bigger than $8$ and that there are only 33 different Ehrhart polynomials\
 corresponding to $\sa_1(P)\in\{1,\dots,35\}\setminus\{33,34\}$.   
 
In dimension 4 we have just made some calculations with the 1561 reflexive  
simplices (cf.~\cite{conrads:reflexive}). Here the Ehrhart polynomials of ''only'' 574 of them have 
roots with real part $-1/2$. Finally we present two 4-dimensional 
reflexive simplices which show that both conditions in Proposition 
\ref{prop:ref_roots} are necessary.  The first simplex is given by the 
inequalities $E_1=\{x\in\R^4 : x_i\geq -1,\,1\leq i\leq 3,\, 
-x_3-2\,x_4\leq,\, 2, \,x_1+x_2+2\,x_3+2\,x_4\leq 1\}$.  With the help 
of the computer program {\tt latte} \cite{latte:program}, which we have used for all our 
calculations, one (the computer) can easily 
determine the Ehrhart polynomial of such a polytope and here we find  
\begin{equation*} 
 \LE(s,E_1) = \frac{27}{2}\,s^4+27\,s^3+21\,y^2+\frac{15}{2}\,y+1. 
\end{equation*}    
Thus we have $\LE(E_1)=70$ and hence $(\LE(E_1)-1-4\,\vol(E_1))^2\geq 
16\,\vol(E_1)$ but $2\,\LE(E_1)>9\vol(E_1)+18$.  
Next let $E_2=\{x\in\R^4 : -x_1\leq 1,\,-x_2\leq 
1,\,-2\,x_1-3\,x_2-4\,x_3\leq 1, \, -4\,x_1-5\,x_2-8\,x_4\leq 1, 
10\,x_1+9\,x_2+4\,x_3+8\,x_4\leq 1\}$. Then  
\begin{equation*} 
 \LE(s,E_2) = \frac{4}{3}\,s^4+\frac{8}{3}\,s^3+\frac{8}{3}\,y^2+\frac{4}{3}\,y+1. 
\end{equation*} 
In this case we have $\LE(E_2)=9$, $2\,\LE(E_2)\leq 9\,\vol(E_2)+18$ but 
$(\LE(E_2)-1-4\,\vol(E_2))^2< 16\,\vol(E_2)$.

\section{\label{sec:3polytopes}$3$-dimensional lattice polytopes}
In this section we will study the roots of Ehrhart polynomials of
$3$-dimensional lattice polytopes. To this end we will distinguish
polytopes with and without interior lattice points. 

\begin{theorem}  Let $\Gamma(3,0)$ be the set of all roots of Ehrhart polynomials of 
$3$-dimensional lattice polytopes $P\in\mathcal{P}^3(0)$, i.e., without interior lattice points. 
\hfill \renewcommand{\labelenumi}{\rm \roman{enumi})} 
\begin{enumerate}
\item $\Gamma(3,0)\cap\R=\{-3,-2\}\cup (-2,1)$. Moreover, $1$ is a
  cluster point and there are  infinitely many roots in  the interval
  $(-2,-1)$. 
\item  $\{a+\mathrm{i}\,b \in \Gamma(3,0) : b\ne 0\}\subset
  W :=\{a+\mathrm{i}\,b : (a+1)^2+b^2\leq 2\text{ and } a\geq -1\}$.
\item On the boundary of the semicircle $W$ lie exactly 33 pairs of
  zeros. $-1\pm \mathrm{i}\,\sqrt{2}$, $-1\pm
  \mathrm{i}/\sqrt{2}$, $-1\pm \mathrm{i}$ and 
   $-1\pm \mathrm{i}/\sqrt{5}$ are the only complex
  roots in $\Gamma(3,0)$ with real part $-1$.
\end{enumerate} 
\label{thm:3dim_without} 
\end{theorem}

For the proof we need the following proposition
\begin{proposition} Let $P\in\mathcal{P}_n$ and let $k\in\N$ be the smallest
  positive integer with $\LE(k\,\inter P)\ne 0$. Then 
\begin{equation*} 
   \LE_{n-1}(P)\leq \frac{n\,k}{2}\,\LE_n(P)=\frac{n\,k}{2}\,\vol(P).
\end{equation*}
\label{prop:general_k}
\end{proposition}  
\begin{proof} Let $P=\{x\in\R^n : a_j\,x\leq b_j,\, 1\leq j\leq m\}$
  be a lattice polytope with facets $F_j$ corresponding the outer normal
  vector $a_j$.  It was already shown by Ehrhart \cite{ehrhart:polynomialII} that 
\begin{equation}
   \LE_{n-1}(P) = \frac{1}{2} \sum_{i=1}^m \frac{\vol_{n-1}(F_i)}{\det(\aff F_i\cap\Z^n)},
\label{eq:2leadterm}
\end{equation}
where $\vol_{n-1}()$ denotes the $(n-1)$-dimensional volume and 
$\det(\aff F_i\cap\Z^n)$ denotes the determinant of the $(n-1)$-dimensional sublattice of $\Z^n$ contained in the affine hull of the facet $F_i$.

Since $P$ is a lattice polytope we can assume $a_j \in\Z^n$, $0\in P$, $b_j\in\N$, and that the
  vectors $a_j$ are primitive, i.e., $\conv\{0,a_j\}\cap\Z^n=
  \{0,a_j\}$. Hence $\det(\aff F_j\cap\Z^n) = \enorm{a_j}$, where
  $\enorm{\cdot}$ denotes the Euclidean norm. By the choice of $k$ we can
  find a $z\in\Z$ such that $(1/k)\,z\in\inter P$ and so we find (cf.~\eqref{eq:2leadterm})
\begin{equation*}
\begin{split}
 \vol(P) & = \frac{1}{n}\sum_{i=1}^m
 \vol_{n-1}(F_i)\,\frac{|a_j\,(1/k)z-b_j|}{\enorm{a_j}} \geq
 \frac{2}{n\,k}\, \frac{1}{2}\sum_{i=1}^m
 \frac{\vol_{n-1}(F_i)}{\enorm{a_j}} \\ &= \frac{2}{n\,k}\,\LE_{n-1}(P)
\end{split} 
\end{equation*}
\end{proof}
We remark that we always have $k\leq n+1$ and thus by Proposition \ref{prop:general_k} 
$\LE_{n-1}(P)\leq\binom{n+1}{2}\vol(P)$ which is a special case of
another series of inequalities proved in \cite{betke-mcMullen:85}. The case $k=1$ and thus $\LE_{n-1}(P)  \leq (n/2) \vol(P)$ was already shown in \cite{wills:81}. 
So with the notation of Proposition \ref{prop:general_k} we have for three-dimensional polytopes $P$ 
\begin{equation}
   1\leq \LE_2(P) \leq k\,\frac{3}{2}\vol(P),
\label{eq:g_2}
\end{equation} 
where the lower bound follows from \eqref{eq:2leadterm} and the fact that for any facet $\vol_{n-1}(F_i)/\det(\aff F_i\cap\Z^n)\geq 1/(n-1)!$.

\begin{proof}[Proof of Theorem \ref{thm:3dim_without}] From 
  \cite[Theorem 1.2, Proposition 4.7]{becketall:roots_ehrhart} it 
  follows that all real roots of Ehrhart polynomials of 3-dimensional
  polytopes are within $[-3,1)$ and in  \cite[Theorem
  1.7]{henkschuermannwills:ehrhart} it was shown that $1$ is cluster
  point of $\Gamma(3,0)$. 
Next we observe that -1 is a root of $\LE(s,P)$ for any
polytope without interior lattice points (cf.~\eqref{eq:reciprocitylaw}). Hence, denoting for
short the coefficients $\LE_i(P)$ by $\lE_i$ we  have
$\lE_3-\lE_2+\lE_1-1=0$ and so may write 
\begin{equation*}
\begin{split}
  \LE(s,P) & = \lE_3\,s^3+\lE_2\,s^2+\lE_1\,s + 1 
             = \lE_3\, (s+1)\left( s^2+\frac{\lE_2-\lE_3}{\lE_3}\,s+\frac{1}{\lE_3}\right).
\end{split} 
\end{equation*}   
For the two remaining roots $-\gamma_{1,2}$ we find 
\begin{equation}
   -\gamma_{1,2} = - \frac{\lE_2-\lE_3}{2\lE_3}\pm
   \sqrt{\left(\frac{\lE_2-\lE_3}{2\,\lE_3}\right)^2-\frac{1}{\lE_3}}.
\label{eq:2-roots}
\end{equation}  
Now we want to show that there are no real roots in $(-3,-2)$. Suppose
$-2$ is another root of $\LE(s,P)$, then for the third root $\gamma$,
say, we get $\gamma=-1/(2\,\lE_3)$. Since $6\,\lE_3$ is an integer we
conclude that $\gamma=-3$ or $\gamma\geq -3/2$. Hence if there is an
Ehrhart-polynomial having a real root in $(-3,-2)$ then we know
$\LE(2\,\inter P)\ne 0$ and so by \eqref{eq:g_2} $\lE_2\leq
3\lE_3$. For given $\lE_3$ the right hand side in \eqref{eq:2-roots}
becomes minimal if $\lE_2$ is as large as possible. Thus
$-\gamma_{1,2}\geq -1\pm\sqrt{1-1/\lE_3} > -2 $. Observe that
$\lE_3\geq 1$ since we have
assumed that all roots are real and $\lE_2\leq
3\lE_3$. 

For i) it remains to show that there are infinitely many real roots in
$(-2,-1)$. To this end we consider for an integer $q$ the 
pyramids $P(q)=\conv\{0,2\,e_1,q\,e_2,\allowbreak 2\,e_1+q\,e_2,
e_3\}$. Then one gets $\LE_3(P(q))=2/3\,q$ and $\LE_2(P(q))=3/2\,q$
which shows by \eqref{eq:2-roots} that for $q$ large $\LE(s,P(q))$ has
a real root in $(-2,-1)$ depending on $q$.

For ii) we assume that the roots $-\gamma_{1,2}$ in \eqref{eq:2-roots}
are complex. Writing 
$-\gamma_{1,2}=a\pm\mathrm{i}\,b$ leads to $ b^2 = 1/\lE_3-a^2$. Since
$1/\lE_3=(1-2\,a)/\lE_2$ we may rewrite this as 
\begin{equation*}
  \left(a+\frac{1}{\lE_2}\right)^2 + b^2 = \left(\frac{1}{\lE_2}\right)^2+ \frac{1}{\lE_2}.
\end{equation*} 
By \eqref{eq:g_2} we know $\lE_2\geq 1$ and it is not hard to see that
all the circles above are contained in the disk given by the largest
one, i.e., we have  $(a+1)^2+b^2\leq 2$. Since we assume that
the roots $-\gamma_{1,2}$ are complex we have $\LE(2\,\inter P)\ne 0$,
because otherwise $-2$ would be a root. Thus from \eqref{eq:g_2} we
conclude $\lE_2\leq 3\,\lE_3$ which is equivalent to  $a=-(\lE_2-\lE_3)/(2\,\lE_3)\geq -1$.
 
Now we come to part iii). Let $\lE_3=\vol(P)=k/6$, $k\in\N$. All
complex roots on the semicircle satisfy $\lE_2=1$ and  
\begin{equation*} 
  0> \left(\frac{\lE_2-\lE_3}{2\,\lE_3}\right)^2-\frac{1}{\lE_3} =
  \left(\frac{3}{k}-\frac{1}{2}\right)^2-\frac{6}{k} =\frac{1}{4\,k^2}\left(k^2-36\,k+36\right).
\end{equation*} 
Hence $k$ is restricted to the integers $k=2,\dots,34$. The
Reeve-simplices $T(k)=\conv\{0,e_1,e_2,(1,1,k)^\intercal\}$ form a
family of simplices whose Ehrhart polynomials 
\begin{equation*}
   \LE(S,T(k))=\frac{k}{6}\,s^3+s^2+\frac{12-k}{6}\,s+1
\end{equation*} 
have these roots. 

Finally we consider the case that the complex roots have real part -1. Then
$\lE_2=3\lE_3$ and the Ehrhart polynomial of such a polytope $P$ is of the type 
\begin{equation*} 
   \LE(s,P) = \lE_3\,s^3+3\,\lE_3\,s^2+(2\,\lE_3+1)\,s +1.
\end{equation*} 
The roots of that polynomial are given by $-1$, $-1\pm\sqrt{1-1/\lE_3}$. 
Again let $\lE_3=k/6$, $k\in \N$. Since $1-1/\lE_3$ has to be negative and since $\lE_3=\lE_2/3\geq 1/3$ we just have to consider the cases $k=2,\dots,5$. Moreover we note that for such a polytope $P$ all roots of $2\,P$ have real part $-1/2$ and so $2P$ has to be a reflexive polytope (cf.~Proposition \ref{prop:roots_reflexive}). Hence all possible candidates  are contained in database of Kreuzer and Skarke of $3$-dimensional reflexive polytopes.

An example for $k=2$ is given by the Reeve-simplex $T(2)$ with Ehrhart polynomial
$(1/3)\,s^3+s^2+5/3\,s+1$ and with complex roots $-1\pm\mathrm{i}\sqrt{2}$. 
For $k=3,4$ we found respectively the simplices $\conv\{0,e_1,e_2,(2,2,3)^\intercal\}$ with  complex roots  $-1\pm \mathrm{i}$ and for $k=4$ the simplex 
$\conv\{0,e_1,e_2,(2,3,4)^\intercal\}$ and complex roots  
$-1\pm 1/\sqrt{2}\,\mathrm{i}$. For $k=5$, i.e., $\lE_3=5/6$, 
there does not exist a simplex with the required Ehrhart polynomial. However, the pyramid over a quadrangle given by $\conv\{0,e_1,2\,e_2, 2\,e_1+e_2,e_3\}$  has the Ehrhart polynomial $5/6\,s^3+5/2\,s^2+16/6\,s+1$ with complex roots $-1\pm \mathrm{i}/\sqrt{5}$.

\end{proof}

Next we come to 3-dimensional polytopes with interior lattice
points. For those lattice polytopes we know by Proposition \ref{prop:general_k} that 
\begin{equation}
     \LE_2(P) \leq \frac{3}{2} \LE_3(P).
\label{eq:interior_points} 
\end{equation}  

First we state some simple properties on the real parts of the roots.

\begin{proposition} Let $P\in\mathcal{P}_n(l)$, $l\geq 1$. \renewcommand{\labelenumi}{\rm \roman{enumi})}  
\begin{enumerate}
\item If 
  all roots of $\LE(s,P)$ are real  then either all roots are contained in
  $(-1,0)$ or one belongs to $(-1,0)$ and the two others are in
  $(0,1)$.
\item If $\LE(s,P)$ has only one real root $\gamma$ then
  $\gamma\in(-1,0)$ and the real parts of the complex roots are contained in   
$(-3/4,1/2)$. 
\end{enumerate} 
\label{prop:realroots}
\end{proposition}
\begin{proof} 
  Let us assume that all roots are real. 
  The point of inflexion of the real polynomial
  $\LE(t,P)$, $t\in\R$, is given by $-\LE_2(P)/(3\,\LE_3(P))$ which by
  \eqref{eq:interior_points} is contained in $[-1/2,0)$. Furthermore
  the derivative of that polynomial at $0$ is given by $\LE_1(P)$ and
  we also know that   $\LE(-1,P)=-l<0$, $\LE(1,P)>0$. Thus, the polynomial has   always
  a real root in $(-1,0)$. If all roots are real then two cases occur.  
  If $G_1(P)\geq 0$ then all of them  
  are in $(-1,0)$ and otherwise one root is contained in
  $(-1,0)$ and the positive roots are strictly less than $1$.

  Now suppose we have one real root $\gamma$ and the two complex roots
  $a\pm\mathrm{i}b$. Since
  $(1/3)(2\,a+\gamma)=-\LE_2(P)/(3\LE_3(P)\in[-1/2,0)$  and
  $\gamma\in(-1,0)$ we must have $-3/4<a<1/2$.
\end{proof}

For the proof of Theorem \ref{thm:roots} we also need the following lemma.

\begin{lemma}  Let $P\in\mathcal{P}_n(l)$, $l\geq 1$. Then \renewcommand{\labelenumi}{\rm \roman{enumi})} 
\begin{enumerate}
 \item $\LE_1(P) \leq
   \LE_2(P)+\LE_3(P)+\frac{2}{3}\leq\frac{5}{2}\LE_3(P)+\frac{2}{3}$,
\item $\LE(-1/(3\,\vol(P),P)\geq 0$,
\end{enumerate}
and both bounds are tight. In particular, equality in ii) is only
attained if $P$ is unimodular equivalent to $S_3(1)$.
\label{lem:dim3_inequalities}
\end{lemma}  
\begin{proof} By \eqref{eq:reciprocitylaw} we have
  $\LE_1(P)=l-\LE_3(P)+\LE_2(P)+1$. By  Theorem
  \ref{thm:non-symmetric} we also know $l\leq 2\,\LE_3(P)/2-1/3$ and
  thus $\LE_1(P)\leq \LE_2(P)+\LE_3(P)+\frac{2}{3}$. The second
  inequality in i) is just a consequence of
  \eqref{eq:interior_points}. 

For the proof of ii) we write for short $\lE_i$ instead of $\LE_i(P)$. On account
of i) we get 
\begin{equation*}
\begin{split}
\LE\left(-\frac{1}{3\vol(P)},P\right) & =-\frac{1}{27(\lE_3)^2}
+\frac{\lE_2}{9(\lE_3)^2}-\frac{\lE_1}{3\lE_3}+1 \\ 
 & \geq -\frac{1}{27(\lE_3)^2}
+\frac{\lE_2}{9(\lE_3)^2}-\frac{\lE_2+\lE_3+2/3}{3\lE_3}+1 \\
& = \frac{1}{27(\lE_3)^2}\left(3\lE_2-1-6\lE_3\right)-\frac{\lE_2}{3\lE_3}+\frac{2}{3}. 
\end{split} 
\end{equation*} 
With $\lE_2\geq 1$ (cf.~\eqref{eq:g_2}) and $\lE_2/(3\lE_3)\leq 1/2$ (cf.~\eqref{eq:interior_points}) we obtain 
\begin{equation*}
  \LE\left(-\frac{1}{3\vol(P)},P\right) \geq \frac{1}{27(\lE_3)^2}
  \left(2-6\lE_3\right)+\frac{1}{6}= \frac{2}{3}\left(\frac{1}{3\,\lE_3}-
      \frac{1}{2}\right)^2 \geq 0.
\end{equation*} 
Now as the inequalities show we have equality in i) if and only if
$\vol(P)=(3\,l+1)6$, i.e., if we have equality in Theorem \ref{thm:non-symmetric}.  In ii)
we have equality if and only if in addition $\vol(P)=2/3$, i.e., $l=1$
and so $P=S_3(1)$. 
\end{proof}

\begin{proof}[Proof of Theorem \ref{thm:roots}] In view of Theorem \ref{thm:3dim_without} and
  Proposition \ref{prop:realroots} it
  remains to show that the norm of each complex root of the Ehrhart
  polynomial of a polytope with interior lattice points is bounded by
  $\sqrt{3}$. Let $-\gamma_1$ be the real root and $a\pm\mathrm{i}\,b$
  be the complex roots with $b\ne 0$. Since $\LE(-1,p)=-l<0$ we get
  from  Lemma \ref{lem:dim3_inequalities} i) that $\gamma_1\geq 1/(3\vol(P))$. On the
  other hand we know that $\gamma_1\cdot(a^2+b^2)=1/\vol(P)$ (cf.~\eqref{eq:coeffroots}) and thus 
  $(a^2+b^2)\leq 3$.

  Among the polytopes $P\in\mathcal{P}_3(l)$, $l\leq 1$, the bound on the
  norm is attained if and only if the polynomial  has two complex roots $a\pm\mathrm{i}\,b$ 
  and one real root $-\gamma_1$ (cf.~Proposition \ref{prop:realroots}). Since  $\gamma_1\cdot(a^2+b^2)=1/\vol(P)$ we get $-\gamma_1=-1/(3\vol(P))$. Thus by Lemma  \ref{lem:dim3_inequalities} ii) we conclude that this is only the
  case for a polytope unimodular equivalent to $S_3(1)$. By Proposition \eqref{prop:non-symmetric}
  we have $\LE(s,S_3(1))=2/3\,s^3+s^2+7/3\,s+1$  with roots $-1/2,-1/2\pm\mathrm{i}\,\sqrt{11}/2$.
\end{proof} 
\vspace{1cm}



\begin{thebibliography}{10}

\bibitem{batyrev:reflexive}
V.V. Batyrev, \emph{Dual polyhedra and mirror symmetry for calabi-yau
  hypersurfaces in toric varieties}, J. Algebr. Geom. \textbf{3} (1994),
  493--535.

\bibitem{becketall:roots_ehrhart}
M.~Beck, J.~{De Loera}, M.~Develin, J.~Pfeifle, and R.P. Stanley,
  \emph{Coefficients and roots of {E}hrhart polynomials}, Contemp. Math.
  \textbf{374} (2005), 15--36.

\bibitem{beckpixton:birkhoff}
M.~Beck and D.~Pixton, \emph{The {E}hrhart polynomial of the {B}irkhoff
  polytope}, Discrete Comput.~Geom. \textbf{30} (2003), no.~4, 623--637.

\bibitem{beckrobins:discrete_continously}
M.~Beck and S.~Robins, \emph{Computing the continuous discretely: Integer-point
  enumeration in polyhedra}, Springer, (to appear), Preprint at
  \url{http://math.sfsu.edu/beck/papers/ccd.html}.

\bibitem{beck_sottile:irrational}
M.~Beck and F.~Sottile, \emph{Irrational proofs for three theorems of
  {S}tanley}, to appear in European J.~Combinatorics, preprint at
  \url{http://front.math.ucdavis.edu/math.CO/0501359}.

\bibitem{betkegritzmann:valuations}
U.~Betke and P.~Gritzmann, \emph{An application of valuation theory to two
  problems of discrete geometry}, Discrete Math. \textbf{58} (1986), 81--85.

\bibitem{betke-mcMullen:85}
U.~Betke and P.~McMullen, \emph{Lattice points in lattice polytopes}, Monatsh.
  Math. \textbf{99} (1985), no.~4, 253--265.

\bibitem{braun:ehrhart_roots}
B.~Braun, \emph{Norm bounds for {E}hrhart polynomials}, preprint at
  \url{http://arxiv.org/abs/math.CO/0602464}.

\bibitem{bump:crosspolytope}
D.~Bump, K.-K. Choi, P.~Kurlberg, and J.~Vaaler, \emph{A local {R}iemann
  hypothesis {I}}, Math. Z. \textbf{233} (2000), no.~1, 1--19.

\bibitem{conrads:reflexive}
H.~Conrads, \emph{Weighted projective spaces and reflexive simplices}, Manuscr.
  Math. \textbf{107} (2002), no.~2, 215--227.

\bibitem{latte:program}
J.~{De Loera}, D.~Haws, R.~Hemmecke, and P.~Huggins, \emph{A user's guide for
  \url{latte} v1.1, software package \url{latte}}, 2004, available at
  \url{http://www.math.ucdavis.edu/~latte}.

\bibitem{diazrobins:latticepolytope}
R.~Diaz and S.~Robins, \emph{The {E}hrhart polynomial of a lattice polytope},
  Ann.~of Math. \textbf{145} (1997), no.~3, 503--518, Erratum in 146:1 (1997),
  237.

\bibitem{ehrhart:polynomial}
E.~Ehrhart, \emph{Sur les poly\`edres rationnels homoth\'etiques \`a n
  dimensions}, C. R. Acad. Sci., Paris, S\'er. A \textbf{254} (1962), 616--618.

\bibitem{ehrhart:polynomialII}
\bysame, \emph{Sur un probl\`eme de g\'eom\'etrie diophantienne lin\'eaire},
  J.~{R}eine {A}ngew.~Math. \textbf{227} (1967), 25--49.

\bibitem{ehrhart:reciprocity}
\bysame, \emph{Sur la loi de r\'eciprocit\'e des poly\`edres rationnels}, C. R.
  Acad. Sci., Paris, S\'er. A \textbf{266} (1968), 695--697.

\bibitem{gruber-lekkerkerker:87}
P.~M. Gruber and C.~G. Lekkerkerker, \emph{Geometry of numbers}, second ed.,
  vol.~37, North-Holland Publishing Co., Amsterdam, 1987.

\bibitem{henkschuermannwills:ehrhart}
M.~Henk, A.~Sch{\"u}rmann, and J.M. Wills, \emph{Ehrhart polynomials and
  successive minima}, to appear in Mathematika, preprint at
  \url{http://front.math.ucdavis.edu/math.MG/0507528}.

\bibitem{hensley:lattice_volume}
D.~Hensley, \emph{Lattice vertex polytopes with interior lattice points}, Pac.
  J. Math. \textbf{105} (1983), 183--191.

\bibitem{hibi:92}
T.~Hibi, \emph{Dual polytopes of rational convex polytopes}, Combinatorica
  \textbf{12} (1992), no.~2, 237--240.

\bibitem{Hibi:lower_bound}
\bysame, \emph{A lower bound theorem for {E}hrhart polynomials of convex
  polytopes}, Adv.~Math. \textbf{105} (1994), no.~2, 162--165.

\bibitem{kirschenhofer_pethoe_tichy:counting_polynomials}
P.~Kirschenhofer, A.~Peth{\"o}, and R.T. Tichy, \emph{On analytical and
  diophantine properties of a family of counting polynomials}, Acta Sci. Math.
  \textbf{65} (1999), 47--59.

\bibitem{kreuzer_skarke:reflexive3}
M.~Kreuzer and H.~Skarke, \emph{Classification of reflexive polyhedra in three
  dimensions}, Adv. Theor. Math. Phys. \textbf{2} (1998), 853--871.

\bibitem{kreuzer_skarke:reflexive4}
\bysame, \emph{Classification of reflexive polyhedra in four dimensions}, Adv.
  Theor. Math. Phys. \textbf{4} (2000), 1209--1230.

\bibitem{lagariasziegler:volume}
J.C. Lagarias and G.M. Ziegler, \emph{Bounds for lattice polytopes containing a
  fixed number of interior points in a sublattice}, Canad.~J.~Math. \textbf{43}
  (1991), 1022--1035.

\bibitem{liu:lattice_face}
F.~Liu, \emph{Ehrhart polynomials of lattice-face polytopes},
  \url{http://arxiv.org/abs/math.CO/0512616}.

\bibitem{liu:cyclic_ehrhart}
\bysame, \emph{Ehrhart polynomials of cyclic polytopes}, J.~Comb.~Theory,
  Ser.~A \textbf{111} (2005), 111--127.

\bibitem{mordell:dedkind}
L.J. Mordell, \emph{Lattice points in tetrahedron and generalized {D}edekind
  sums}, J.~Indian Math.~Soc. (N.S.) \textbf{15} (1951), 41--46.

\bibitem{mustatepayne:bettyehrhart}
M.~Mustata and S.~Payne, \emph{Ehrhart polynomials and stringy {B}etti
  numbers}, \url{http://arxiv.org/abs/math.AG/0505054}.

\bibitem{nill:thesis}
B.~Nill, \emph{Gorenstein toric {F}ano varieties}, Ph.D. thesis, University
  T{\"u}bingen, 2005,
  http://w210.ub.uni-tuebingen.de/dbt/volltexte/2005/1888/pdf/nill.pdf.

\bibitem{pick:lattice}
G.A. Pick, \emph{Geometrisches zur {Z}ahlenlehre}, Sitzungsber.~Lotus Prag
  \textbf{19} (1899), 311--319.

\bibitem{pikhurko:lattice_volume}
O.~Pikhurko, \emph{Lattice points in lattice polytopes}, Mathematika
  \textbf{48} (2001), no.~1-2, 15--24.

\bibitem{pommersheim:tetrahedra}
J.E. Pommersheim, \emph{Toric varieties, lattice points and {D}edekind sums},
  Math.~Ann. \textbf{295} (1993), no.~1, 1--24.

\bibitem{Poonen_Rodriguez-Villegas:reflexive}
B.~Poonen and F.~Rodriguez-Villegas, \emph{Lattice polygons and the number 12},
  Am. Math. Mon. \textbf{107} (2000), no.~3, 238--250.

\bibitem{rodriguez:zeros}
F.~Rodriguez-Villegas, \emph{On the zeros of certain polynomials}, Proc. Amer.
  Math. Soc. \textbf{130} (2002), 2251--2254.

\bibitem{stanley:nonnegative}
R.P. Stanley, \emph{Decompositions of rational convex polytopes}, Ann.~Discrete
  Math. \textbf{6} (1980), 333--342.

\bibitem{stanley:symmetric}
\bysame, \emph{On the number of faces of centrally-symmetric simplicial
  polytoeps}, Graphs and Combinatorics \textbf{3} (1987), 55--66.

\bibitem{stanley:monoton}
\bysame, \emph{A monotonicty property of $h$-vectors and $h^\star$-vectors},
  European J. Combinatorics \textbf{14} (1993), no.~3, 251--258.

\bibitem{wills:81}
J.M. Wills, \emph{On an analog to {M}inkowski's lattice point theorem}, The
  geometric vein, Springer, New York, 1981, pp.~285--288.

\bibitem{wills:90c}
\bysame, \emph{Minkowski's successive minima and the zeros of a
  convexity-function}, Monatsh. Math. \textbf{109} (1990), no.~2, 157--164.

\bibitem{zpw:82}
J.~Zaks, M.~A. Perles, and J.M. Wills, \emph{On lattice polytopes having
  interior lattice points}, Elem. Math. \textbf{37} (1982), no.~2, 44--46.

\end{thebibliography}
\def\cprime{$'$} \def\cprime{$'$} \def\cprime{$'$} \def\cprime{$'$}
  \def\cprime{$'$} \def\cprime{$'$} \def\cprime{$'$}
\providecommand{\bysame}{\leavevmode\hbox to3em{\hrulefill}\thinspace}
\providecommand{\MR}{\relax\ifhmode\unskip\space\fi MR }
\providecommand{\MRhref}[2]{%
  \href{http://www.ams.org/mathscinet-getitem?mr=#1}{#2}
}
\providecommand{\href}[2]{#2}

\end{document}